\documentclass[12pt,eqno]{article}
\usepackage{amsmath,amsthm,amscd,amssymb}
\usepackage{latexsym}
\usepackage{extarrows}
\usepackage{color}
\newtheorem{lemma}{Lemma}[section]
\newtheorem{remark}{Remark}[section]
\newcommand{\bremark}{\begin{remark} \em}
\newcommand{\eremark}{\end{remark} }
\newtheorem{theorem}{Theorem}[section]

\newtheorem{definition}[theorem]{Definition}

\begin{document}
\parindent 14pt
\renewcommand{\theequation}{\thesection.\arabic{equation}}
\renewcommand{\baselinestretch}{1.15}
\renewcommand{\arraystretch}{1.1}
\def\disp{\displaystyle}

\title{\large Weak entropy solution for a Keller-Segel type fluid model   \footnotetext{{\small E-mails: chen@math.uni-mannheim.de, fhuang@amt.ac.cn,  liulingjun@amss.ac.cn }\\ }}
\author{\vspace{-0.4em}{\small Li Chen$^1$, Feimin Huang$^2$, Lingjun Liu$^2$ }\\
\vspace{-0.4em}{\small 1. School of Business Informatics and Mathematics, University of Mannheim, }\\
\vspace{-0.4em}{\small
68159 Mannheim, Germany}\\
\vspace{-0.4em}{\small 2. Academy of Mathematics and Systems Science, Academia Sinica,}\\
{\small Beijing 100190, P.R.China.}\\
 \date{}}
\maketitle
\vspace{-3em}
\abstract{\small In this paper, we consider a Keller-Segel type fluid model, which is a kind of Euler-Poisson system with a self-gravitational force. We show that similar to the parabolic case, there is a critical mass $8\pi$ such that if the initial total mass $M$ is supercritical, i.e., $M> 8\pi$, then any weak entropy solution with the same mass $M$ must blow up in finite time.
The a priori estimates of weak entropy solutions for critical mass $M=8\pi$ and subcritical mass $M<8\pi$ are also obtained.}

\


\section{ Introduction}\label{section1}

\setcounter{section}{1}

\setcounter{equation}{0}

Chemotaxis stands
	for the movement of organism in response to a chemical stimulus, this phenomenon can
	be observed from microscopic bacteria to the largest mammals. 
 The Keller-Segel system proposed in \cite{KS,KS1} is the basic model to describe chemotaxis. Since the beginning of this century, the analysis on critical mass has attacted more and more experts to join the research of investigating the wellposedness on different variant Keller-Segel type systems, \cite{AT, BCL, BCM, CCLW, CLW12, CW14, CS, JL, LS06, N,Sugiyama06, SS, TW, TW1, WCH}, to name a few. The critical mass of the simplified parabolic elliptic  two-dimensional Keller-Segel system was studied in \cite{BDP}, i.e.,
\begin{equation}\label{eq01}
\left\{
\begin{array}{ll}
\partial_t\rho=\Delta \rho-\nabla \cdot (\rho\nabla \Phi), \quad x\in\mathbb{R}^2,\quad t>0,\\
-\Delta\Phi=\rho,
\end{array}
\right.
\end{equation}
where $\rho(x, t)$ represents the cell density, and $\Phi(x,t)$ is the concentration of chemo-attractant which induces a drift force.
It was shown in \cite{DP} that there exists a critical mass $8\pi$ such that any solution in the supercritical mass case blows up in finite time.
After the blow-up time, $\delta$-measure is formed.
The global existence and large time behavior of weak solutions for both the critical and subcritical cases were also obtained, for example in \cite{BDP,BCM}. Recently, a microscopic understanding of the Keller-Segel system is given through the mean field limit of a system of interacting particle system in \cite{AG,HL}. In the microscopic level, the case with heavy particle dynamics for the mass accumulation \cite{AG}, which corresponds to the blow up in the macroscopic level, is still unsolved.

In this paper, we consider a related model in the form of Euler-Poisson system, i.e.,
\begin{equation}\label{eq1}
\left\{
\begin{array}{ll}
\partial_t\rho+\nabla\cdot m=0, \quad x\in\mathbb{R}^2,\quad t>0,\\
\partial_tm+\nabla\cdot(\frac{m}{\rho}\otimes m)+\nabla \rho=\rho\nabla\Phi-m,\\
-\Delta\Phi=\rho,
\end{array}
\right.
\end{equation}
where $m$ is the momentum of the fluids, and denote $u=\frac{m}{\rho}$ for $\rho>0$ as usual.
	
The microscopic understanding of system \eqref{eq1} is the following. We start with an interacting $N$ particle system based on Newton's second law and the dynamic driven by mean field force
	\begin{eqnarray}
		\left\{
\begin{array}{ll}
		\dfrac{d}{dt}X_i(t)=V_i(t),\quad 1\leq i\leq N,\label{ParticalSystem}\\
		\dfrac{d}{dt}V_i(t)=\dfrac{1}{N}\sum^N_{j=1,j\neq i}\nabla \Psi(X_i(t)-X_j(t))-\frac{1}{\tau}V_i(t),\nonumber
	\end{array}
	\right.
	\end{eqnarray}
	where the mean field interaction force is given by the Coulomb potential $\nabla \Psi(x)=\nabla (-\frac{1}{2\pi}\log |x|)=-\frac{x}{2\pi|x|^2}$ in 2-D.
	Through the characteristic method, the empirical measure $f^N(t,x,v)$ given by $\frac{1}{N}\sum^N_{i=1}\delta_{X_i(t),V_i(t)}(x,v)$, where $(X_i(t),V_i(t))$ is the solution of \eqref{ParticalSystem}, formally satisfies the following Vlasov Poisson equation in the sense of distribution,
	\begin{eqnarray}
		\partial_t f^N+v\cdot\nabla_x f^N+\nabla_v\cdot ((\nabla \Psi * \rho^N-v)f^N)=0,\label{VlasovType}
	\end{eqnarray}
	where $\rho^N(t,x)=\int_{\mathbb R^2} f^N(t,x,v)dv$.
	Then after taking $N\rightarrow \infty$ formally, the so called mean field limit, the above equation \eqref{VlasovType} is named as the mean field equation. We refer to \cite{LP,JW} for rigorous proof of this limit.
	If furthermore we take the local Maxwellian Ansatz in \eqref{VlasovType},
	$$
	f(t,x,v)=\rho(t,x)\frac{1}{(2\pi)^\frac{d}{2}}e^{-\frac{|v-u(t,x)|^2}{2}}, \quad d=2,
	$$ 
	with $\rho(t,x)=\int_{\mathbb R^2} f(t,x,v)dv$ and $u(t,x)=\int_{\mathbb R^2} vf(t,x,v)dv$, we arrive at system \eqref{eq1} with $m=u{\rho}$. 
	If one takes the relaxation limit in \eqref{eq1} formally, one can obtain $m=\rho\nabla\Phi-\nabla\rho$. Hence together with the conservation of mass one obtains exactly the parabolic-elliptic Keller-Segel equation \eqref{eq01}.

 It is convenient to write the third equation in \eqref{eq1} in the compact form \begin{equation}\label{Phi}\Phi(x,t)= \mathbb{G}*\rho(x,t)=-\frac{1}{2\pi}\int_{\mathbb{R}^2}\log|x-y|\rho(y,t)dy,\end{equation} where $ \mathbb{G}(x)=-\frac{1}{2\pi}\log|x|$ is the Green's function for $-\Delta$ in $\mathbb{R}^2$. We assume that the initial data $\rho_0(x):=\rho(x,0),  m_0(x):=m(x,0)$ 
satisfies
 \begin{equation}\label{eq2}
	\begin{aligned}
		\rho_0(1+|x|^2)+\frac{|m_0|^2}{\rho_0},~ \rho_0\log\rho_0\in L^1(\mathbb{R}^2).
	\end{aligned}
	\end{equation}
	The conserved total mass is
\begin{equation}\label{eq3}
\begin{aligned}
M:=\int_{\mathbb{R}^2}\rho(x,t)dx=\int_{\mathbb{R}^2}\rho_0(x)dx .
\end{aligned}
\end{equation}

	First we recall a logarithmic Hardy-Littlewood-Sobolev inequality. 

\vskip 0.1in

\begin{lemma}\cite{Be} \label{lhs}
	 Suppose that $f$ be a non-negative function in $L^1(\mathbb{R}^2)$ and $\int_{\mathbb{R}^2}fdx=M$ such that $f\log f$ and $f\log(1+|x|^2)$ belong to $L^1(\mathbb{R}^2)$, then
	 \begin{equation}\label{eq3b0}
	 \begin{aligned}
	 \int_{\mathbb{R}^2}f\log fdx+\frac{2}{M}\iint\limits_{\mathbb{R}^2\times\mathbb{R}^2}f(x)f(y)\log|x-y|dxdy\geq-\mathfrak{C}(M),
	 \end{aligned}
	 \end{equation}
	 with $\mathfrak{C}(M):=M(1+\log \pi-\log M)$.
\end{lemma}

It should be noted from \eqref{eq1}$_3$ that the initial data $\Phi_0(x)$ of $\Phi(x,t)$ is completely determined by the initial density $\rho_0(x)$. In fact, multiplying \eqref{eq1}$_3$ for $t=0$ by $\Phi_0$  yields that 
\begin{equation}\label{eq11b11}
\begin{aligned}
\int_{\mathbb{R}^2}|\nabla\Phi_0|^2dx=\int_{\mathbb{R}^2}\rho_0\Phi_0 dx
=-\frac{1}{2\pi}\iint\limits_{\mathbb{R}^2\times\mathbb{R}^2}\rho_0(x)\rho_0(y)\log|x-y|dxdy. 
\end{aligned}
\end{equation}
Then using Lemma \ref{lhs}, we have 
\begin{equation}\label{eq1c1b11}
	\begin{aligned}
		\int_{\mathbb{R}^2}\rho_0\log\rho_0 dx-\frac{4\pi}{M}\int_{\mathbb{R}^2}|\nabla\Phi_0|^2dx\geq-\mathfrak{C}(M),
	\end{aligned}
	\end{equation}
which deduces that 
	\begin{equation}\label{eq1cc1b11}
		\begin{aligned}
			\int_{\mathbb{R}^2}|\nabla\Phi_0|^2dx\leq\frac{M}{4\pi}\big(\mathfrak{C}(M)+\int_{\mathbb{R}^2}\rho_0\log\rho_0 dx\big). 
		\end{aligned}
		\end{equation}	
	\begin{remark} 
The initial data \eqref{eq2} implies that the initial total energy $E_0$ is finite, where 
		\begin{equation}\label{eq02}
	\begin{aligned}
	E_0:=3\int_{\mathbb{R}^2}\frac{|m_0|^2}{\rho_0}+2\rho_0\log\rho_0-|\nabla\Phi_0|^2dx+ \int_{\mathbb{R}^2}\rho_0|x|^2+2x\cdot m_0dx.
	\end{aligned}
	\end{equation}
	 	
	\end{remark}



Note that the system \eqref{eq1} admits a mechanical entropy-entropy flux pair which is convex and reads 
\begin{equation}\label{eq110}
\begin{aligned}
\eta(\rho,m)=\frac{|m|^2}{\rho}+2\rho \log\rho, \qquad
q(\rho,m)=\frac{|m|^2m}{\rho^2}+2m\log\rho+2m.
\end{aligned}
\end{equation}
We define the weak entropy solution of the Euler-Poisson equations \eqref{eq1} with the initial values \eqref{eq2} as follows.
\vskip 0.1in

\begin{definition}[\bf Weak entropy solution]\label{def1}
	A triplet of measurable functions $(\rho, m, \Phi)$  are called weak entropy solution of the Cauchy problem \eqref{eq1} and \eqref{eq2} provided that
	\begin{equation}\label{eq30}
	\left\{
	\begin{array}{ll}
	\displaystyle\int_0^\infty\int_{\mathbb{R}^2}(\rho\psi_t+m\cdot\nabla\psi )dxdt+\int_{\mathbb{R}^2}\rho_0(x)\psi(x,0)dx=0,  \\[4mm]
	\displaystyle\int_0^\infty\int_{\mathbb{R}^2}\big[m\psi_t+\frac{m}{\rho}\otimes m):\nabla\psi+\rho\nabla\psi- m\psi+\rho\nabla\Phi\psi\big] dxdt+\int_{\mathbb{R}^2}m_0(x)\psi(x,0)dx=0,\\[4mm]
	\displaystyle\int_0^\infty\int_{\mathbb{R}^2}\nabla\Phi\cdot\nabla\psi+\rho\psi dxdt=0,\\[4mm]
	\displaystyle\int_0^\infty\int_{\mathbb{R}^2}(\eta\phi_t+q\cdot\nabla\phi-\frac{|m|^2}{\rho}\phi+m\cdot\nabla\Phi\phi)dxdt+\int_{\mathbb{R}^2}\eta(x,0)\phi(x,0)dx\geq0,
	\end{array}
	\right.
	\end{equation}
	holds for any test function $\phi(>0),\psi\in C^\infty_0({\mathbb{R}^2}\times[0,+\infty))$.

\end{definition}
Since the system \eqref{eq1} is a multi-dimensional hyperbolic system of conservation laws with self-gravitational force, singularity like shock might be formed in finite time. The global existence of weak entropy solution in the sense of Definition \ref{def1} is a challenging problem. For the works on Euler-Poisson system with repulsive force, there are few results about $N=2$. Nevertheless, we try to extend the works of \cite{BDP} in this paper to the  system \eqref{eq1} and give at the same time a priori estimates for the weak entropy solution. We assume that the weak entropy solutions $(\rho,m,\Phi)(x,t)$ of \eqref{eq1} and \eqref{eq2} satisfy
\begin{equation}\label{eq2a}
\begin{aligned}
\rho(1+|x|^2)+\frac{|m|^2}{\rho}, \rho\log\rho\in C\big(0,T; L^1(\mathbb{R}^2)\big),
\end{aligned}
\end{equation}
for any time $T>0$.
Then we have

\vskip 0.1in

\begin{theorem}[\bf Supercritical mass]\label{theorem1}
	Assume that the total mass $M> 8\pi$ and the weak entropy solution $(\rho, m, \Phi)$ of \eqref{eq1} and \eqref{eq2} satisfies \eqref{eq3},  \eqref{eq2a} and the following boundary condition, 
	\begin{equation}\label{boundary}
	(|x|^2\rho, m,\Phi)(x,t)\to (0,0,0),~~\mbox{as}~~ |x|\to \infty.
	\end{equation}
If there exists a constant $C>0$ such that 
	\begin{equation}\label{eq6}
	\begin{aligned}
	\int_{\mathbb{R}^2}\rho(x,t) \log\rho(x,t) dx\leq Ct^\alpha, \quad 0<\alpha<1,
	\end{aligned}
	\end{equation}
then the  weak entropy solution   $(\rho, m, \Phi)$ blows up in finite time, i.e., 
$$\exists \ T^*>0, \quad s.t \quad \int_{\mathbb{R}^2}|x|^2\rho(x,T^*)dx=0.$$


\end{theorem}
\begin{remark}
	The growth condition \eqref{eq6} for density $\rho$ is mild since it is uniformly bounded with respect to $t$ for subcritical mass below. In the parabolic case, which is a time relaxation version of this fluid model, the blow up appears at a time when both the physical entropy and the potential energy are going to infinity. Therefore intuitively if the physical entropy grows in time more than linearly, i.e., the assumption \eqref{eq6} does not hold, then the vanishing of the second moment will not happen in finite time.
\end{remark}

\begin{remark}
Note that the weak entropy solution belongs to a natural functional space, which allows the discontinuity of  the density, and the total mass is conserved, Theorem \ref{theorem1} may hint the formation of Delta-measure in finite time $(t=T^*)$ like the Keller-Segel system \eqref{eq01}. 
\end{remark}

\begin{theorem}[\bf  Critical and subcritical mass]\label{theorem2}
	 Assume that $(\rho, m, \Phi)$ is any weak entropy solution  of \eqref{eq1} and \eqref{eq2} satisfying \eqref{eq3},  \eqref{eq2a} and \eqref{boundary}. Then 
	 
	 (i) If $M=8\pi$ (critical mass), it holds that for any $T>0$, 
	 \begin{equation}\label{eq6b1}
	 \begin{aligned}
	 \int_{\mathbb{R}^2}\frac{1}{2} |x|^2\rho(x,T)dx+\int_{\mathbb{R}^2}\frac{|m(x,T)|^2}{\rho(x,T)}dx+\iint\limits_{\mathbb{R}^2\times[0,T]}\frac{2|m|^2}{\rho}dxdt\leq  C_1(M)+E_0;
	 \end{aligned}
	 \end{equation}
	 
	 (ii) If $M<8\pi$ (subcritical mass), it holds that for any $T>0$,
	 \begin{equation}\label{eq6d1}
	 \begin{aligned}
	 \int_{\mathbb{R}^2} \frac{|m(x,T)|^2}{\rho(x,T)}dx+\iint\limits_{\mathbb{R}^2\times[0,T]}\frac{2|m|^2}{\rho}dxdt
	 \leq C_2(M)(1+\log(10+T)+E_0),
	 \end{aligned}
	 \end{equation}
	 \begin{equation}\label{eq32p}
	 \begin{aligned} C_3(M)-E_0-C_4(M)\log(10+T)\leq\int_{\mathbb{R}^2}\rho(x,T)\log\rho(x,T) dx\leq C_5(M), 
	 \end{aligned}
	 \end{equation}
	 and
	 \begin{equation}\label{eq6c1}
	 \begin{aligned}
	 &\int_{\mathbb{R}^2} |x|^2\rho(x,T)dx
	 \leq C_6(M)(1+T+E_0),
	 \end{aligned}
	 \end{equation}
	 where $E_0$ is the initial total energy and $C_i(M), i=1,\cdots,6$ are positive constants which only depend on  $M$, but are independent of $T$. 
	
\end{theorem}
\begin{remark}
	Theorem \ref{theorem2} can be regarded as a priori estimates of weak entropy solutions, which might be useful for the existence and stability problems in the cases of critical and subcritical mass. 
\end{remark}

The rest of the paper is arranged as follows. Some useful energy inequalities are derived in section \ref{section2} from  an intrinsic property of weak  solutions and the entropy inequality 
$\eqref{eq30}_3$. Then the blow-up phenomenon for supercritical mass is investigated in section \ref{section3} and the a priori estimates for critical and subcritical mass are obtained in section \ref{section4}. 

\vskip 0.1in










\section{ Energy inequalities}\label{section2}

\setcounter{section}{2}

\setcounter{equation}{0}

\noindent 
Motivated by \cite{BDP}, we have an intrinsic property of weak solutions, i.e., 
\begin{lemma}\label{lemma2}
 Assume that $(\rho, m, \Phi)$ is any weak entropy solution  of \eqref{eq1} and \eqref{eq2}  satisfying \eqref{eq3},  \eqref{eq2a} and \eqref{boundary}, 
%
	it holds that
	\begin{equation}\label{eq4}
	\begin{aligned}
	\int_{\mathbb{R}^2}(|x|^2\rho_t+2x\cdot m_t)dx=4M(1-\frac{M}{8\pi})+\int_{\mathbb{R}^2}\frac{2|m|^2}{\rho}dx.
	\end{aligned}
	\end{equation}
	

\end{lemma}
\vskip 0.1in

\begin{proof}

From the system \eqref{eq1}$_2$, we can obtain that
 \begin{equation}\label{eq23}
\begin{aligned}
m=\rho\nabla\Phi-\nabla \rho-\partial_tm-\nabla\cdot(\frac{m}{\rho} \otimes m).
\end{aligned}
\end{equation}
Taking $\nabla\cdot$\eqref{eq23} and inserting the result into \eqref{eq1}$_1$, we have
 \begin{equation}\label{eq24}
\begin{aligned}
\partial_t\rho=\Delta\rho-\nabla\cdot(\rho\nabla\Phi)+\nabla\cdot[\partial_tm+\nabla\cdot(\frac{m}{\rho} \otimes m)].
\end{aligned}
\end{equation}
We choose a smooth function $\varphi_\varepsilon(|x|)$ with compact support that grows nicely to $|x|^2$ as $\varepsilon\rightarrow0$, that is
\begin{equation}\label{eq28}
\begin{aligned}
\lim\limits_{\varepsilon\rightarrow0}\varphi_\varepsilon(|x|)=|x|^2,\quad\quad \lim\limits_{\varepsilon\rightarrow0}\nabla\varphi_\varepsilon(|x|)=2x,\quad \quad\lim\limits_{\varepsilon\rightarrow0}\Delta\varphi_\varepsilon(|x|)=4 .
\end{aligned}
\end{equation}

Since $\Delta\varphi_\varepsilon$ is bounded and $\nabla\varphi_\varepsilon(x)$ is Lipschitz continuous, we can use Definition \ref{def1} to get
\begin{equation}\label{eq7}
\begin{aligned}
&\frac{d}{dt}\int_{\mathbb{R}^2}\varphi_\varepsilon\rho dx\\
=&\int_{\mathbb{R}^2}\Delta  \varphi_\varepsilon\rho dx+\int_{\mathbb{R}^2}\nabla\varphi_\varepsilon\cdot(\rho\nabla\Phi)dx-\int_{\mathbb{R}^2}\nabla\varphi_\varepsilon\cdot[m_t+\nabla\cdot(\frac{m}{\rho} \otimes m)]dx.
\end{aligned}
\end{equation}
Using the symmetry of \eqref{Phi} yields that 
\begin{equation}\label{eq25}
\begin{aligned}
\int_{\mathbb{R}^2}\nabla\varphi_\varepsilon\cdot(\rho\nabla\Phi)dx & =-\frac{1}{2\pi}\int\limits_{\mathbb{R}^2\times\mathbb{R}^2}\nabla\varphi_\varepsilon(x)\cdot\nabla \log|x-y|\rho(x,t)\rho(y,t)dxdy \\ & =-\frac{1}{4\pi}\iint\limits_{\mathbb{R}^2\times\mathbb{R}^2}\frac{(\nabla\varphi_\varepsilon(x)-\nabla\varphi_\varepsilon(y))\cdot(x-y)}{|x-y|^2}\rho(x,t)\rho(y,t)dxdy.
\end{aligned}
\end{equation}
Thus, we deduce that
\begin{equation}\label{eq8}
\begin{aligned}
\frac{d}{dt}&\int_{\mathbb{R}^2}(\varphi_\varepsilon\rho+\nabla\varphi_\varepsilon\cdot m )dx\\
&=\int_{\mathbb{R}^2}\Delta \varphi_\varepsilon\rho dx-\frac{1}{4\pi}\iint\limits_{\mathbb{R}^2\times\mathbb{R}^2}\frac{(\nabla\varphi_\varepsilon(x)-\nabla\varphi_\varepsilon(y))\cdot(x-y)}{|x-y|^2}\rho(x,t)\rho(y,t)dxdy \\ & \qquad \quad \qquad \quad-\int_{\mathbb{R}^2}\nabla\varphi_\varepsilon\cdot\nabla\cdot(\frac{m}{\rho} \otimes m)dx.
\end{aligned}
\end{equation}
Sending $\varepsilon\rightarrow0$ gives that

\begin{equation}\label{eq29}
\begin{aligned}
\frac{d}{dt}\int_{\mathbb{R}^2}&(|x|^2\rho+2x\cdot m) dx+\int_{\mathbb{R}^2}2x\cdot\nabla\cdot(\frac{m}{\rho} \otimes m)dx \\
 &=\int_{\mathbb{R}^2}4\rho dx-\frac{1}{4\pi}\iint\limits_{\mathbb{R}^2\times\mathbb{R}^2}2\rho(x,t)\rho(y,t)dxdy \\
& =4M(1-\frac{M}{8\pi}).
\end{aligned}
\end{equation}

On the other hand, a direct computation gives that 
\begin{equation}\label{eq40}
\begin{aligned}
\int_{\mathbb{R}^2}2x\cdot\nabla\cdot(\frac{m}{\rho} \otimes m)dx=-2\int_{\mathbb{R}^2}\frac{|m|^2}{\rho}dx.
\end{aligned}
\end{equation}
Then we have
\begin{equation}\label{eq31}
\begin{aligned}
\frac{d}{dt}\int_{\mathbb{R}^2}(|x|^2\rho+2x\cdot m) dx=4M(1-\frac{M}{8\pi})+\int_{\mathbb{R}^2}\frac{2|m|^2}{\rho}dx,
\end{aligned}
\end{equation}
which completes the proof of Lemma \ref{lemma2}.

\end{proof}


From Lemma \ref{lemma2} and the entropy inequality $\eqref{eq30}_3$, we have 
\begin{lemma}\label{lemma3}
	Assume that $(\rho, m, \Phi)$ is any weak entropy solution  of \eqref{eq1} and \eqref{eq2}  satisfying \eqref{eq3},  \eqref{eq2a} and \eqref{boundary}, it holds that  for any $T>0$, 
\begin{equation}\label{ineq1}
\begin{aligned}
\displaystyle\int_{\mathbb{R}^2}\frac{1}{2}|x|^2\rho(x,T) dx\leq&4M(1-\frac{M}{8\pi})T+\int_{\mathbb{R}^2}\frac{2|m(x,T)|^2}{\rho(x,T)}dx\\
& +\iint\limits_{\mathbb{R}^2\times[0,T]}\frac{2|m|^2}{\rho}dxdt
+\int_{\mathbb{R}^2}(\rho_0|x|^2+2x\cdot m_0)dx,
\end{aligned}
\end{equation}
and 
\begin{equation}\label{ineq2}
\begin{aligned}
\int_{\mathbb{R}^2}\frac{|m(x,T)|^2}{\rho(x,T)}dx&+\iint\limits_{\mathbb{R}^2\times[0,T]}\frac{2|m|^2}{\rho}dxdt+2(1-\frac{M}{8\pi})\int_{\mathbb{R}^2}\rho(x,T)\log\rho(x,T) dx\\
\quad \ \leq& \frac{M}{4\pi}\mathfrak{C}(M)+\int_{\mathbb{R}^2}\frac{|m_0|^2}{\rho_0}+2\rho_0 \log\rho_0-|\nabla\Phi_0|^2dx, 
\end{aligned}
\end{equation}
where $\mathfrak{C}(M)$ is a constant defined in Lemma \ref{lhs}. 
\end{lemma}

\begin{remark}
	Estimate \eqref{ineq2} shows that in the subcritical case that the kinetic energy and the physical entropy $\int \rho\log\rho$ is uniform bounded in time, while in the supercritical case the kinetic energy is dominated by the physical entropy.
\end{remark}

\begin{proof}
Integrate \eqref{eq4} over $[0,T]$, we obtain that 
\begin{equation}\label{eq32}
\begin{aligned}
&\int_{\mathbb{R}^2}(|x|^2\rho(x,T)+2x\cdot m(x,T))dx\\
\ \ =&4M(1-\frac{M}{8\pi})T+\iint\limits_{\mathbb{R}^2\times[0,T]}\frac{2|m|^2}{\rho}dxdt+\int_{\mathbb{R}^2}(\rho_0|x|^2+2x\cdot m_0)dx,
\end{aligned}
\end{equation}
which directly gives \eqref{ineq1}  from the following Cauchy inequality  
\begin{equation}\label{eq3a2}
\begin{aligned}
\int_{\mathbb{R}^2}|2x\cdot m| dx\leq\int_{\mathbb{R}^2}(\frac{1}{2}|x|^2\rho+\frac{2|m|^2}{\rho})dx.
\end{aligned}
\end{equation}

It remains to show \eqref{ineq2}. Note that the inequality \eqref{eq30}$_3$ is equivalent to
\begin{equation}\label{eq11}
\begin{aligned}
\partial_t\eta+\nabla\cdot q\leq2m\cdot\nabla\Phi-\frac{2|m|^2}{\rho}
\end{aligned}
\end{equation}
in the sense of distribution. That is,
\begin{equation}\label{eq111}
\begin{aligned}
\int_{\mathbb{R}^2}(\partial_t\eta&+\nabla\cdot q)dx=\frac{d}{dt}\int_{\mathbb{R}^2}(\frac{|m|^2}{\rho}+2\rho \log\rho) dx\\
&\leq\int_{\mathbb{R}^2}(2m\cdot\nabla\Phi-\frac{2|m|^2}{\rho})dx=\frac{d}{dt}\int_{\mathbb{R}^2}|\nabla\Phi|^2dx-\int_{\mathbb{R}^2}\frac{2|m|^2}{\rho}dx.
\end{aligned}
\end{equation}
Integrating \eqref{eq111} over $[0,T]$ gives that 
\begin{equation}\label{eq1111}
\begin{aligned}
\int_{\mathbb{R}^2}&\frac{|m(x,T)|^2}{\rho(x,T)}+\iint\limits_{\mathbb{R}^2\times[0,T]}
\frac{2|m|^2}{\rho}dxdt+\int_{\mathbb{R}^2}(2\rho(x,T) \log\rho(x,T)- |\nabla\Phi(x,T)|^2)dx\\
\  \ \leq&\int_{\mathbb{R}^2}(\frac{|m_0|^2}{\rho_0}+2\rho_0 \log\rho_0-|\nabla\Phi_0|^2)dx.
\end{aligned}
\end{equation}

On the other hand, multiplying \eqref{eq1}$_3$ by $\Phi$  yields that 
\begin{equation}\label{eq11b1b1}
\begin{aligned}
\int_{\mathbb{R}^2}|\nabla\Phi(x,T)|^2dx=\int_{\mathbb{R}^2}\rho(x,T)\Phi(x,T) dx
=-\frac{1}{2\pi}\iint\limits_{\mathbb{R}^2\times\mathbb{R}^2}\rho(x,T)\rho(y,T)\log|x-y|dxdy, 
\end{aligned}
\end{equation}
which, together with Lemma \ref{lhs}, gives that 
\begin{equation}\label{eq11c11}
\begin{aligned}
&\int_{\mathbb{R}^2}2\rho(x,T) \log\rho(x,T)-|\nabla\Phi(x,T)|^2dx\\
=&\int_{\mathbb{R}^2}2\rho(x,T)\log\rho(x,T) dx+\frac{1}{2\pi}\iint\limits_{\mathbb{R}^2\times\mathbb{R}^2}\rho(x,T)\rho(y,T)\log|x-y|dxdy\\
=&\frac{M}{4\pi}\big(\int_{\mathbb{R}^2}\rho(x,T) \log\rho(x,T)dx+\frac{2}{M}\iint\limits_{\mathbb{R}^2\times\mathbb{R}^2}\rho(x,T)\rho(y,T)\log|x-y|dxdy\big)\\
&\qquad \quad+2(1-\frac{M}{8\pi})\int_{\mathbb{R}^2}\rho(x,T) \log\rho (x,T)dx\\
 \geq& 2(1-\frac{M}{8\pi})\int_{\mathbb{R}^2}\rho(x,T) \log\rho(x,T) dx-\frac{M}{4\pi}\mathfrak{C}(M).
\end{aligned}
\end{equation}
The inequality \eqref{ineq2} immediately holds from \eqref{eq1111} and \eqref{eq11c11}. Thus the proof of Lemma \ref{lemma3} is completed. 
%
\end{proof}

\section{ Supercritical mass}\label{section3}
\setcounter{section}{3}
\setcounter{equation}{0}

This section is devoted to the supercritical mass, i.e.,  $M>8\pi$. In this case, $\frac{M}{4\pi}-2>0$. 

\begin{proof}[\bf Proof of Theorem \ref{theorem1}]
From the assumption \eqref{eq6} and \eqref{ineq2}, we have
\begin{equation}\label{eq1o111}
\begin{aligned}
&\int_{\mathbb{R}^2}\frac{|m(x,T)|^2}{\rho(x,T)}dx+\iint\limits_{\mathbb{R}^2\times[0,T]}\frac{2|m|^2}{\rho}dxdt\\
\leq& \int_{\mathbb{R}^2}(\frac{|m_0|^2}{\rho_0}+2\rho_0 \log\rho_0-|\nabla\Phi_0|^2)dx+(\frac{M}{4\pi}-2)\int_{\mathbb{R}^2}\rho(x,T) \log\rho(x,T) dx+\frac{M}{4\pi}\mathfrak{C}(M)\\
\leq& \frac{M}{4\pi}\mathfrak{C}(M)+C(\frac{M}{4\pi}-2)T^\alpha+ \int_{\mathbb{R}^2}\frac{|m_0|^2}{\rho_0}+2\rho_0 \log\rho_0-|\nabla\Phi_0|^2dx,\\
\end{aligned}
\end{equation}
for $0<\alpha<1$. 
	From the Lemma \ref{lemma3}, let \eqref{ineq1}$+$\eqref{ineq2}$\times 2$, we get 
\begin{equation}\label{ineq11}
	\begin{aligned}
	\displaystyle\int_{\mathbb{R}^2}\frac{1}{2}&|x|^2\rho(x,T) dx+\iint\limits_{\mathbb{R}^2\times[0,T]}\frac{2|m|^2}{\rho}dxdt+4(1-\frac{M}{8\pi})\int_{\mathbb{R}^2}\rho(x,T)\log\rho(x,T) dx\\
	&\leq4M(1-\frac{M}{8\pi})T+\frac{M}{2\pi}\mathfrak{C}(M)
	+\int_{\mathbb{R}^2}\rho_0|x|^2+2x\cdot m_0dx\\
	&\qquad \qquad +2\int_{\mathbb{R}^2}\frac{|m_0|^2}{\rho_0}+2\rho_0 \log\rho_0-|\nabla\Phi_0|^2dx.
	\end{aligned}
	\end{equation}
Inserting \eqref{eq1o111} into \eqref{ineq11} implies that
\begin{equation}\label{eq3s2}
\begin{aligned}
\int_{\mathbb{R}^2}\frac{1}{2}|x|^2\rho(x,T)dx\leq4M(1-\frac{M}{8\pi})T+C(\frac{M}{8\pi}-1)T^\alpha+\tilde{C}_1(M)+E_0,
\end{aligned}
\end{equation}
where $E_0$ is the initial total energy and $\tilde{C}_1(M)>0$ is a constant which only depends on $M$.

%
Note that $0<\alpha<1$ and $M>8\pi$, there exists 
a finite time $T^*>0$ such that the right hand side of \eqref{eq3s2} is negative as $T>T^*$, which contradicts with the fact that the left hand side of  \eqref{eq3s2} is always non-negative for any time $T>0$. Therefore there is no global existence of weak entropy solution satisfying \eqref{eq2a} and \eqref{eq6} in the case of supercritical mass. And the entropy solution $(\rho, m, \Phi)$ blows up in finite time.
	\end{proof}
\section{ Critical and subcritical mass }\label{section4}

\setcounter{section}{4}

\setcounter{equation}{0}

\noindent Since the time-asmptotical behavior of \eqref{eq1} (the relaxation limit) is formally dominated by the Keller-Segel system \eqref{eq01},  we expect that the system \eqref{eq1} admits a global weak entropy solution for both critical and subcritical mass. 
In \cite{BDP}, the global existence of the Keller-Segel system \eqref{eq01} was shown by the celebrated Aubin-Lions Lemma, in which the derivative estimates concerning with $\rho|\nabla \log\rho|^2$ play essential role. However, 
it is difficult to achieve such derivative estimates for the system \eqref{eq1} since it is essentially hyperbolic and shock may occur. The method of Aubin-Lions Lemma may not be available to study the global existence of weak entropy solution for the system \eqref{eq1}. Nevertheless, we can obtain some a priori  estimates for both critical and subcritical mass through Lemma \ref{lemma3}, and  expect that these estimates might be useful for the global existence of weak entropy solution of \eqref{eq1}.

\begin{proof}[\bf Proof of Theorem \ref{theorem2}]

\

\begin{itemize} 
\item {\bf Critical mass, i.e. $M=8\pi$. }

Note that  $M=8\pi$, \eqref{eq6b1} immediately holds from \eqref{ineq1} and \eqref{ineq2}. 

%

\item 
{\bf Subcritical mass, i.e. $M<8\pi$. } 

We first show the lower bound of $\int_{\mathbb{R}^2}\rho(x,T) \log\rho(x,T)dx$.
Let $\mu(x)=\frac{1}{\pi(10+T)}e^{-\frac{|x|^2}{10+T}}$, then $\int_{\mathbb{R}^2}\mu(x) dx=1$ and $\int_{\mathbb{R}^2}d\nu(x)=1$ with $d\nu(x)=\mu(x)dx$. Note that $\rho\log\rho$ is convex, the Jensen's inequality gives that 
\begin{equation}\label{eq3a6}
\begin{aligned}
 &\int_{\mathbb{R}^2}\rho(x,T)\log\rho(x,T) dx\\
=&\int_{\mathbb{R}^2}\rho(x,T) \log\rho(x,T)+\rho(x,T)\log e^\frac{|x|^2}{10+T} dx-\frac{1}{10+T}\int_{\mathbb{R}^2}|x|^2\rho(x,T) dx\\
=&\int_{\mathbb{R}^2}\rho(x,T) \log(\frac{\rho(x,T)}{\mu})dx-\int_{\mathbb{R}^2}\rho \log[\pi(10+T)]dx-\frac{1}{10+T}\int_{\mathbb{R}^2}|x|^2\rho(x,T) dx\\
=&\int_{\mathbb{R}^2}\frac{\rho(x,T)}{\mu} \log(\frac{\rho(x,T)}{\mu}) d\nu-M\log[\pi(10+T)]-\frac{1}{10+T}\int_{\mathbb{R}^2}|x|^2\rho (x,T)dx\\
\geq&(\int_{\mathbb{R}^2}\frac{\rho(x,T)}{\mu} d\nu)\log(\int_{\mathbb{R}^2}\frac{\rho(x,T)}{\mu}d\nu)-M\log[\pi(10+T)]-\frac{1}{10+T}\int_{\mathbb{R}^2}|x|^2\rho(x,T) dx\\
=&M\log M-M\log[\pi(10+T)]-\frac{1}{10+T}\int_{\mathbb{R}^2}|x|^2\rho(x,T) dx.
\end{aligned}
\end{equation}
Note that $M<8\pi$, then $1-\frac{M}{8\pi}>0$. Inserting \eqref{eq3a6} into \eqref{ineq2} infers that 
\begin{equation}\label{eq1d111}
\begin{aligned}
&\int_{\mathbb{R}^2}\frac{ |m(x,T)|^2}{\rho(x,T)}dx+\iint\limits_{\mathbb{R}^2\times[0,T]}\frac{2 |m|^2}{\rho}dxdt\\
\le& \tilde{C}_2(M)+M\log(10+T)+\frac{1}{10+T}\int_{\mathbb{R}^2}|x|^2\rho(x,T) dx,
\end{aligned}
\end{equation}
 which, together with \eqref{ineq1}, gives \eqref{eq6c1}, which claims $\int_{\mathbb{R}^2}|x|^2\rho(x,T)dx$ is bounded by $1+T$, where $\tilde{C}_2(M)>0$ only depends on $M$ and the initial data. Again using \eqref{eq1d111} yields \eqref{eq6d1}. 
%
%

On the other hand,  the upper bound of 
$\int_{\mathbb{R}^2}\rho\log\rho dx$ in the inequality \eqref{eq32p} can be directly obtained by \eqref{ineq2} due to the fact that $1-\frac{M}{8\pi}>0$. The lower bound is achieved from \eqref{eq3a6} and \eqref{eq6c1}. Therefore the proof of Theorem \ref{theorem2} is completed. 

%
\end{itemize}
\end{proof}

\noindent {\bf Acknowledgments:}
F. Huang's research is supported in part by the
National Natural Science Foundation of China No. 11371349 and 11688101.

\end{document}